\documentclass[12pt]{amsart} 

\usepackage{amscd,amssymb,mathtools,enumerate}
\usepackage[plainpages,urlcolor=blue,backref]{hyperref}

\theoremstyle{plain}
\newtheorem{thm}[subsection]{Theorem}
\newtheorem{lem}[subsection]{Lemma}
\newtheorem{prop}[subsection]{Proposition}
\newtheorem{cor}[subsection]{Corollary}
\newtheorem{conj}[subsection]{Conjecture}

\theoremstyle{definition}

\newtheorem{ex}[subsection]{Example}

\numberwithin{equation}{section}
\newcommand{\bb}{\mathbb}

\newcommand{\coker}{\text{Coker\ }}

\newcommand{\defect}{\mathrm{def}}

\newcommand{\mc}{\mathcal}

\newcommand{\p}{\partial}

\newcommand{\st}{\mathrm{\ s.t.\ }}

\newcommand{\Sym}{\mathrm{Sym}}

\newcommand{\wh}{\widehat}

\newcommand{\xto}{\xlongrightarrow}

\begin{document}
\date{}

  \title[New proofs for technical results]{New proofs for technical results in ``Infinitesimal invariants of mixed Hodge structures'' (arXiv:2406.17118v1) }

   \author[ZHENJIAN WANG]{ZHENJIAN WANG}
   \address{Zhenjian Wang, Hefei National Laboratory, Hefei 230088, China}
   \email{wzhj@ustc.edu.cn}



\begin{abstract} 
Cubic forms $C$ are constructed in the work of R. Aguilar, M. Green and P. Griffiths to establish the generic global Torelli theorem for Fano-K3 pairs $(X,Y)$, where $X: F=0$ is a cubic threefold in $\bb{P}^4$ and $Y\in|-K_X|$ is an anticanonical smooth section of $X$ defined by a quadratic form $Q$. In this article, we prove the following two results, which were previously verified with the computer aid of Macaulay2: for a generic pair $(X,Y)$, (i) the cubic form $C$ is smooth; (2) $(J_{F,3}:Q)=0$, and thereby give a precise meaning of the word ``generic" in this context. 
\end{abstract}


\maketitle


\section{Statements of main results}
In the recent paper \cite{AGG24}, the authors study infinitesimal invariants for pairs $(X, Y)$, where $X$ is a cubic threefold and $Y$ is a smooth anticanonical K3 section. 
They establish a generic global Torelli theorem for such pairs, with a certain cubic form playing a central role in the proof.

Following the notation of \cite[Section 4]{AGG24}, let $V$ be a five‑dimensional complex vector space and let $X \subset \bb{P}V^* \cong \bb{P}^4$ be a cubic threefold defined by ${F=0}$. Take $Y = \{F=0, Q=0\}$ and assume both $X$ and $Y$ are smooth.
Denote by $J_F$ the ideal in $S^\bullet V =\Sym^\bullet V$ generated by the first partial derivatives of $F$, and set $J_{F,k}=J_F\cap S^kV$. The cubic form $C$ is constructed in \cite[Section 3.2]{AGG24} via the composition
\[
H^1(\Omega_X^2)\otimes H^0(-K_X)\to H^1(T_X)
\]
and the derivative of the period map 
\[
H^1(T_X)\to S^2H^1(\Omega_X^2)^*
\]
(cf. \cite[Section 5.5]{CMP17}), followed by
\[
H^1(\Omega_X^2)^*\otimes S^2H^1(\Omega^2_X)^*\to S^3H^1(\Omega_X^2)^*,
\]
yielding a map $H^0(-K_X)\to S^3H^1(\Omega^2_X)^*$. Equivalently, the cubic form $C$ can be characterized by the relation $\bb{C}C = (J_{F,5}:Q)^\perp$ (see \cite[Lemma 4.8]{AGG24}).

The following technical statements are proved in \cite{AGG24} with the computer aid of Macaulay2.

\begin{thm}[{\cite[Theorem 1.1, Theorem 4.7]{AGG24}}]\label{tec1}
For $X$ a cubic threefold and $Y$ an anticanonical smooth section,
both of them generic, the cubic $C$ is smooth.
\end{thm}
\begin{prop}[{\cite[Proposition 4.17]{AGG24}}]\label{tec2}
For generic $F\in S^3V,Q\in S^2V$, we have $(J_{F,3}:Q)=0$, where
\[
(J_{F,3}:Q)=\{a\in V: aQ\in J_{F,3})\}.
\]
\end{prop}

The proofs of these results rely heavily on explicit computer algebra, which makes it difficult to extend the discussion in \cite{AGG24} to more general pairs $(F,Q)$ of higher degree or in more variables. This limits the applicability of the general theory developed there. Moreover, the exact meaning of the term ``generic'' in the statements remains unclear.

Inspired by the work in \cite{Wan15}, which gives a precise interpretation of ``generic'' in \cite[Section 4, (b)]{CG80}, we prove analogous results that serve the same purpose for Theorem \ref{tec1} and Proposition \ref{tec2}.

Let us fix notation, the duality pairing between $S^3V$ and $S^3V^*$ is denoted simply  by $\langle\ ,\ \rangle: S^3V\times S^3V^*\to\bb{C}$. For any subset $E\subset S^3V$, set
\[
E^\perp=\{a\in S^3V^*:\langle b,a\rangle=0,\forall b\in E\rangle\}.
\]
Denote by $\mc{U}\subset S^3V$ the set of smooth cubic forms $F$ such that $\left(J_{F,3}\right)^\perp\subset S^3V^*$ contains a smooth cubic.

\begin{thm}\label{thm:main1}
With the notation above, $\mc{U}$ is a nonempty Zariski open subset of $S^3V$. Moreover, for every $F\in\mc{U}$, there exists $Q\in S^2V$, such that the surface $Y=\{F=Q=0\}\subset\bb{P}V^*$ is smooth and the associated cubic $C$ is smooth.
\end{thm}
\begin{thm}\label{thm:main2}
With the notation above, for every smooth $F\in S^3V$, there exists $Q\in S^2V$, such that  $Y=\{F=Q=0\}\subset\bb{P}V^*$ is smooth and $(J_{F,3}:Q)=0$.
\end{thm}

These results yield a more precise version of the technical statements in \cite{AGG24}.
\begin{cor}
With the notation above, given an arbitrary $F\in\mc{U}$, for generic $Q\in S^2V$, 
\begin{enumerate}[(i)]
\item the pair $(F,Q)$ is a Fano-K3 pair;
\item the cubic form $C$ associated with $(F,Q)$ is smooth;
\item $(J_{F,3}:Q)=0$.
\end{enumerate}
\end{cor}

It would be interesting to give a more concrete description of the set $\mc{U}$ and  to relate it to polynomials of Sebastiani-Thom type, as in \cite{UY09, Wan15}. 

\bigskip
We would like to thank Professor A. Dimca for bringing the work \cite{AGG24} to our attention. This work is supported by NSFC-12301100.

\section{Defects of linear systems and strong Lefschetz property}
Let $V$ be a complex vector space of dimension $n+1\geq3$. Up to identification $\bb{P}V^*\cong\bb{P}^n$, the homogeneous coordinate ring $S^\bullet V$ of $\bb{P}V^*$ is identified with $S=\bb{C}[x_0,\ldots, x_n]$ with natural grading $S=\oplus_{d=0}^\infty S_d$; and identify $S^\bullet V^*$ with $\bb{C}[y_0,\ldots, y_n]$.
The perfect pairing $\langle\ ,\ \rangle: S^kV\times S^kV^*\to\bb{C}$ is the same as the polar paring
\[
(F,G)\mapsto G\left(\frac{\p}{\p x_0},\ldots,\frac{\p}{\p x_n}\right)F.
\]
Given a homogeneous polynomial $F\in S_d$, the Jacobian ideal of $F$, denoted by $J_F$, is the graded ideal in $S$ generated by the first partial derivatives $\frac{\p}{\p x_i}, i=0,\ldots, n$, and the Milnor algebra is $M(F)=S/J_F$. 

\subsection{Defects of linear systems} We are going to determine the dimensions $\dim M(F)_k$ for $F\in S_d$. It is clear that $\dim M(F)_k=\dim S_k=\binom{n+k}{k}$ for $k\leq d-2$. If the hypersurface $\{F=0\}\subset\bb{P}^n$ is smooth, then the dimensions $\dim M(F)_k$ are given explicitly in terms of the Hilbert-Poincar\'e series
\[
\text{HP}\left(M(F),t\right)=\sum_{k\geq0}\dim M(F)_kt^k=\frac{(1-t^{d-1})^{n+1}}{(1-t)^{n+1}};
\]
c.f. \cite[Proposition 7.22]{Dim87}. In particular, it follows that $M(F)_k = 0$ for $k > T=(n+1)(d-2)$ and $\dim M(F)_k = \dim M(f_s)_{T-k}$ for $0 \leq k \leq T$. On the other hand, when $\{F=0\}$ is a hypersurface having only isolated singularities, denote by $\Sigma_F$ its singular locus subscheme,
then $\dim M(F)_k$ is closely related to the defect $\defect S_k(\Sigma_F)$ of linear system $S_k(\Sigma_F)$,
\[
\defect_k\, \Sigma_F=\dim H^0(\Sigma_F,\mc{O}_{\Sigma_F})-\dim\frac{S_k}{(\wh{J_F})_k},
\]
where $\wh{J_F}$ is the saturation ideal of $J_F$,
\[
\wh{J_F}=\{s\in S:\forall i=0,\ldots, n,\exists m_i\geq0 \st x_i^{m_i}\in J_F\}.
\]

From Theorem 1 in \cite{Dim13}, we have the following useful formula
\begin{lem}\label{lem:defect} Let $F\in S_d$ be such that $\{F=0\}\subset\bb{P}^n$ has only isolated singularities with singular locus subscheme $\Sigma_F$. Then
\[
\dim M(F)_{T-k} = \dim M(F_s)_k + \defect_k\, \Sigma_F
\]
for $0 \leq k \leq nd - 2n - 1$, where $T = T(n, d) = (n + 1)(d - 2)$, and $F_s\in S_d$ is a homogeneous polynomial such that $\{F_s=0\}$ is a smooth hypersurface in $\mathbb{P}^n$.
\end{lem}

When $\Sigma_F=\{a_1,\ldots, a_p\}$ is a reduced scheme consisting of $p$ points, consider the evaluation map
\[
\theta_k:\ S_k\to\bb{C}^p,\qquad H\mapsto(H(a_1),\ldots, H(a_p)),
\]
then $\defect_k\, \Sigma_F=\dim(\coker\theta_k)$; c.f. \cite[Chapter 6, Exercise 4.4 (i)]{Dim92}.

\begin{ex}\label{ex:defect}
Let
\[
Q(x_0,\ldots, x_n)=\begin{dcases}
\sum\limits_{0\leq i<j<k\leq n} x_ix_jx_k &\text{if }d=3,\\
\sum\limits_{0\leq i<j\leq n}(x_i^{d-2}x_j^2+x_i^2x_j^{d-3}) &\text{if }d\geq4.
\end{dcases}
\]
Then $J_Q$ is contained in the ideal that is generated by all the products $x_ix_j$ for all $0\leq i<j\leq n$. Moreover, in the proof of Proposition 4.3 in \cite{AI14}, it is the shown that $J_{Q,T-1}$ (recall $T=(n+1)(d-2)$) is equal to $\mc{W}$, where $\mc{W}$ is the subspace of $S_{T-1}$ spanned by all monomials in $x_0, \dots, x_n$ of degree $T-1$ other than $x_0^{T-1},\ldots, x_n^{T-1}$. Therefore, the singular point set of $\{Q=0\}$, which is clearly the zero locus of $J_{Q,T-1}$, consists
of the $n+1$ points $p_0,\ldots, p_n$ in $\bb{P}^n$ given by the classes of the canonical basis $e_0,\ldots, e_n$ of the
vector space $\bb{C}^{n+1}$. 

One sees easily that all $p_i$ are nodes (alias, $A_1$ singularities), hence $\{Q=0\}$ is a nodal hypersurface with reduced singular locus $\Sigma_Q=\{p_0,\ldots, p_n\}$.
Furthermore, for the set $\Sigma_Q$, the above map $\theta_k$ is surjective for all $k\geq1$. It follows that $\defect_k\,\Sigma_Q=0$ for all $k$. Thus, we have
\[
\dim M(Q)_k = \dim M(F_s)_k.
\]
for $k\leq T-1$ by Lemma \ref{lem:defect}.
\end{ex}

\subsection{Leftschetz property of the Milnor algebras}
It is a well known conjecture that the Milnor algebra $M(F)$ for a smooth homogeneous polynomial $F$ has the strong Lefschetz
property; see \cite{DI25} for a brief overview.

\begin{conj}\label{conj:Lef}
For any smooth $F\in S_d$, any integer $k\in[0,T/2)$ ($T=(n+1)(d-2)$), and any generic linear form $\ell\in S_1$, the induced multiplication map
\[
\ell^{T-2k}:\ M(F)_k\to M(F)_{T-k}
\]
is an isomorphism.
\end{conj}

To prove our main theorem, we will need the following remarkable result which is proved in \cite[Theorem C]{BFP23}. 
\begin{lem}\label{lem:Lef}
For $n=4$, if $\{F=0\}$ is a smooth cubic threefold, then the Minor algebra $M(F)$ has the strong Lefschetz
property.
\end{lem}

\section{Proofs of main results}
Let $n=4$ and $V$ be a 5-dimensional vector space. Recall that the set $\mc{U}$ is the set of cubic forms $F\in S^3V\cong S_3$ satisfying the following properties:
\begin{enumerate}[(i)]
\item $F$ is smooth, i.e., the hypersurface $X=\{F=0\}\subset\bb{P}V^*$ is smooth;
\item there exists a smooth cubic $G\in S^3V^*$ such that $\langle J_{F,3},G\rangle=0$, i.e., $G\in\left(J_{F,3}\right)^\perp$.
\end{enumerate}

\subsection{Proof of Theorem \ref{thm:main1}} 
We first prove the non-emptiness of $\mc{U}$. Following the idea in the proof of \cite[Proposition 4.3]{AI14}\cite[Theorem 1.3]{DIN25}, we prove that $\mc{U}$ is nonempty by constructing some explicit hypersurfaces with isolated singularities.

For the polynomial $Q=\sum_{0\leq i<j<k\leq 4}x_ix_jx_k$, we have proved that $\dim M(Q)_3=\dim M(F_s)_3=10$ in Example \ref{ex:defect}. Hence $\dim\left(J_{Q,3}\right)^\perp=10$. In addition, it is obvious that the Fermat cubic $\sum_{i=0}^4 y_i^3$ is contained in $\left(J_{Q,3}\right)^\perp$.

Since the set of smooth polynomials is dense in $S^3$, we may choose a sequence of smooth cubics $F_i, i=1,2,\ldots$ converging to $Q$ in $S^3V$. Since $\dim \left(J_{F_i,3}\right)^\perp=10$ for all $i$, passing to a subsequence if necessary, the vector subspaces $\left(J_{F_i,3}\right)^\perp$ of $S^3V^*$ converges to a 10-dimension subspace $E$ of $S^3V^*$. Since $\dim\left(J_{Q,3}\right)^\perp=10$ as well, it follows that $E=\left(J_{Q,3}\right)^\perp$, and hence $E$ contains a smooth polynomial $\sum_{i=0}^4 y_i^3$. Therefore, for $i$ sufficiently large, $\dim \left(J_{F_i,3}\right)^\perp$ also contains a smooth polynomial, since the set of smooth polynomials form an open dense subset; then this $F_i$ belongs to $\mc{U}$. It follows that $\mc{U}$ is a nonempty Zariski open subset of $S^3V$, since $\mc{U}$ is clearly a Zariski open subset. 

Now given $F\in\mc{U}$, and a smooth $G\in\left(J_{F,3}\right)^\perp$. The vector space
\[
G^\perp=\{b\in S^3V: \langle b,G\rangle=0\}
\]
is a hyperplane in $S^3V$ containing $J_{F,3}$, and thus gives a hyperplane of $S^3V/J_{F,3}$. By Macaulay's theorem \cite[Theorem 7.4.1]{CMP17}, we have a nondegenerate pairing $S^2V/J_{F,2}\times S^3V/J_{F,3}\to S^5V/J_{F,5}\cong\bb{C}$ induced by multiplication, and thus the annihilator of the hyperplane induced by $G^\perp$ in $S^2V/J_{F,2}$ is one-dimensional which is denoted by $\bb{C}Q'\mod J_{F,2}$ for $Q'\in S^2V$. It is clear that $(J_{F,5}:Q')=G^{\perp}$.

If $\{F=Q'=0\}\subset\bb{P}V^*$ is smooth, then by \cite[Lemma 4.8]{AGG24}, the cubic form associated to $(F,Q')$ is $(J_{F,5}:Q')^\perp$, which is equal to $G$, hence smooth. Otherwise, it suffices to find a $q\in J_{F,2}$ such that $Q=Q'+q$ satisfies
\begin{enumerate}[(i)]
\item $\{F=Q=0\}\subset\bb{P}V^*$ is smooth;
\item the cubic form associated to $(F,Q)$ is smooth.
\end{enumerate}
Indeed, the linear space spanned by $Q'$ and $J_{F,2}$ induces a base point free linear system on $X=\{F=0\}$ since $F$ is smooth and $J_{F,2}$ is already base point free. It follows that for a general element $Q=Q'+q$ in $Q'+J_{F,3}$, the surface $\{F=Q=0\}$ is smooth by Bertini's theorem \cite{Jou83}. Furthermore, since $q\in J_{F,2}$, we have
\[
(J_{F,5}:Q)=(J_{F,5}:Q'),
\]
and hence by \cite[Lemma 4.8]{AGG24}, the cubic form associated to $(F,Q)$ is $G$ which is smooth. The proof of Theorem \ref{thm:main1} is complete.

\subsection{Proof of Theorem \ref{thm:main2}}
Recall that $X=\{F=0\}$ is a smooth cubic 3-fold. By Lemma \ref{lem:Lef} has the strong Lefschetz property, namely, for a general linear form $\ell\in V$, the multiplications
\[
\ell^3: V\to S^4V/J_{F,4}\qquad\text{and}\qquad \ell:S^2V/J_{F,2}\to S^3V/J_{F,3}
\]
are isomorphisms.

The surface $\{F=Q=0\}$ is smooth for a general $Q\in S^2V$, as it follows from the proof of Theorem \ref{thm:main1}. We claim $(J_{F,3}:Q)=0$ for a general $Q$, i.e., 
\[
V\xto{Q} S^3V/J_{F,3}
\]
is injective for a general $Q$. In fact, by the strong Lefschetz property of the Milnor algebra $S^\bullet V/J_F$, we have that 
\[
V\xto{\ell^2} S^3V/J_{F,3}
\]
is injective for a general $\ell\in V$; hence the injectivity still holds for a general $Q\in S^2V$ by the openness of injective homomorphisms. The proof of Theorem \ref{thm:main2} is complete.


\end{document}